\documentstyle[12pt]{article}
\def\pd#1#2{{\partial#1\over\partial#2}}
\def\pdd#1#2#3{{\partial^2#1\over{\partial#2\partial#3}}}
\def\dint{\int\!\!\!\!\!\int}
\def\E{\mbox{\rm E}}
\def\P{\mbox{\rm P}}
\def\var{\mbox{\rm var}}
\def\cov{\mbox{\rm cov}}
\def\X{{\cal X}}
\def\A{{\cal A}}
\def\B{{\cal B}}
\def\norm#1{\left|\!\left|#1\right|\!\right|}
\def\bds{\begin{displaystyle}}
\def\eds{\end{displaystyle}}
\def\ben{\begin{enumerate}}
\def\een{\end{enumerate}}
\newtheorem{theorem}{Theorem}

\newtheorem{lemma}[theorem]{Lemma}
\begin{document}
\title{\bf  On Solutions of First Order Stochastic Partial Differential Equations}
\author{Hamza K. and Klebaner F.C.\thanks{School of Mathematical Sciences, Monash University, Clayton,
Victoria, AUSTRALIA 3800. Research supported by the Australian
Research Council.
2005. AMS Subject Classification: 60H15.}\\
Monash University}
\date{}
\maketitle

\subsection*{Abstract}
This note is concerned with an important for modelling question of existence of solutions of stochastic partial differential equations as
proper stochastic processes, rather than processes in the generalized sense.
We consider a first order stochastic partial differential equations  of the form $\pd Ut = DW$,
and $\pd Ut-\pd Ux= DW$, where $D$ is a differential operator and $W(t,x)$ is a continuous but
 non-differentiable function (field).
 We give a necessary and sufficient condition for stochastic  equations
to have solutions as functions.
The result is then applied to the  equation for a yield curve. Proofs are based on   probability
arguments.

\section{Introduction}

Stochastic differential equations are often obtained from ordinary ones by introduction of
noise, which is taken to be a generalized function (Schwartz distribution).
Following Rozanov 1995, the  noise in  SPDE is defined as follows.
Let $W(t,x)$ be a continuous  function (field) and $D$   a differential operator.
If $W$ is not differentiable, then $DW$   does not exist in
the usual sense of a function, but  can be  considered as a generalized function. It is defined  by the following action on any  test function $\varphi (t,x)
\in C_K^{\infty}((0,\infty)\times (0,\infty))$, the space of infinitely differentiable compactly supported functions,
$$\left<DW,\varphi\right>=\left<W,D^\star
\varphi\right>=\dint W(t,x) D^\star \varphi (t,x) dtdx,$$
where $D^\star$ is the
adjoint operator of $D$. The adjoint operator is defined by the identity
$\left<Df,\varphi\right>=\left<f,D^\star\varphi\right>$ that
holds for all
smooth functions $f$ and   test functions $\varphi$. For details see Rozanov 1995 p. 99-103, or
 Lang 1993.

In applications the noise is typically random.
White noise, for example, is a generalized derivative  of Brownian motion
$\dot B(t)=\frac{dB(t)}{dt}$, and
the noise  considered in second order partial differential equations
is the second   derivative of Brownian sheet $B(t,x)$,
$\dot B(t,x)=\pdd Btx(t,x)$, (Walsh
1984, Carmona and Nualart 1988, Freidlin 1988.)
First we recall a definition of
 Brownian sheet $B(t,x)$ on positive quadrant of the plane $R_+^2$, which is
 the basic model for White noise field with its various modifications.
In our applications $t$ is time and is non-negative, as well as the state variable $x$. This is not really a
restriction,   if needed  the Brownian sheet can be taken on the whole plane $R^2$.
 A Gaussian random measure $\B$ on $R^2$ is defined
 by the
following properties. For any Borel set $A\subset R_+^2$,
$\B(A)$ is a Normal random variable with zero mean  and  variance given by the area of $A$. For non overlapping
$A_1,A_2$,  $\B(A_1)$ and $\B(A_2)$ are independent and
$\B(A_1\cup  A_2)=\B(A_1)+\B(A_2)$. Put $B(t,x)=\B([0,t]\times ([0,x])$.  It is known that for almost
all realizations, $B(t,x)$ is a continuous  but nowhere differentiable function of $t,x$. For details see e.g.
Walsh 1984.

 In the next section we give results for solutions of first order partial differential equations
 to be functions when the noise is  a generalized function, and
 then use probabilistic arguments to derive conditions for the case when  noise is obtained from a stochastic
 process, such as a white noise. In Section 3 we apply the
 result to the equation  of a yield curve. Section 4 contains the proofs. Since our focus is on modelling
 with SPDE's it is of prime concern that solutions will be proper functions. The classical existence and uniqueness results
 apply in the space of generalized functions. Perhaps a classical approach is
 to use existence and uniqueness in the space of generalized functions and then apply some regularization
 results. Our approach is more direct, and
 utilizes the probabilistic
 nature of noise in the equations.


\section{First Order SPDE's}

Let $D$ be  a first order differential operator
\begin{equation}
D\varphi(t,x) = a(t,x)\pd\varphi t(t,x)+b(t,x)\pd\varphi
x(t,x)+c(t,x)\varphi(t,x),
\label{gdo}
\end{equation}
where   $a$ and $b$ are smooth
functions from $C^1$, and $c$ is continuous. The adjoint operator
$D^\star$  is easily found to be
$$D^\star\varphi = -\pd{(a\varphi)}t-\pd{(b\varphi)}x+c\varphi.$$
It turns out that only
a particular form of the differential operator $D$
yields solutions as functions of the SPDE (\ref{genPDE}).
\begin{theorem}
Let $D$ is given by (\ref{gdo}) and $U$ be a solution to the pde
\begin{equation}
\pd Ut = DW,
\label{genPDE}
\end{equation}
in the sense that
$-\left<U,\pd\varphi t\right> = \left<W,D^\star\varphi\right>.$ for any test function $\varphi$.
Then it holds
\begin{eqnarray}
\lefteqn{\int_0^x[U(t,y)-U(0,y)]dy}\nonumber\\
&&= \int_0^t\left[b(s,x)W(s,x)-b(s,0)W(s,0)\right]ds +
\int_0^x\left[a(t,y)W(t,y)-a(0,y)W(0,y)\right]dy\nonumber\\
&&+\int_0^x\left[\int_0^t\left(\pd at+\pd bx-c\right)(s,y)W(s,y)ds\right]dy.
\label{solgenPDE}
\end{eqnarray}
Thus a solution of (\ref{genPDE}) as a function exists if and only if
$\bds\int_0^tb(s,x)W(s,x)ds\eds$ is
differentiable in
$x$. In particular, when $b$ is identically zero, a solution as a function  exists.
\label{genPDEth}
\end{theorem}

Consider now the equation
\begin{equation}
\pd rt-\pd rx = DW.
\label{gspde}
\end{equation}
We say that $r(t,x)$ is a
solution of (\ref{gspde}) if (\ref{defgspde}) holds for any test
function $\varphi$
\begin{equation}
\dint r(t,x)\left(\pd\varphi t-\pd\varphi x\right)(t,x)dtdx
= \dint W(t,x)\left(\pd{(a\varphi)}t
+\pd{(b\varphi)}x-c\varphi\right)(t,x)dtdx.
\label{defgspde}
\end{equation}

\begin{theorem}
Suppose that the noise process is such that for any $x>0$, $W(0,x)=0$.
Then the equation
\begin{equation}
\pd rt(t,x)-\pd rx(t,x) =  a(t,x)\left(\pd Wt(t,x)-\pd
Wx(t,x)\right)+c(t,x)W(t,x)
\label{pspde}
\end{equation}
(this is equation (\ref{gspde}) with $b(t,x)=-a(t,x)$) with the initial condition $r(0,x) = r_0(x)$
and  $r(t,x)$ is continuous at $t=0$, for all $x$, has for unique solution the
function
\begin{equation}
r(t,x)  = a(t,x)W(t,x)  + \int_0^tW(s,x)\left[\pd ax(s,x)-\pd
at(s,x)+c(s,x)\right]ds + r_0(t+x).
\label{gsol}
\end{equation}
\label{pspdeth}
\end{theorem}


Next (motivated by applications) we take the noise  $W(t,x)$ as the Brownian sheet
\begin{equation}
W(t,x)= B(t,t+x):=B^x(t).
\label{modWN}
\end{equation}
\begin{theorem}
A solution to equation (\ref{gspde}), where $D$ is given by (\ref{gdo}) and $W(t,x)$ by (\ref{modWN}), as a function exists
  if and only if
\begin{equation}
a(t,x)=-b(t,x).
\label{cond}
\end{equation}
Moreover, in this case, the general solution (\ref{gsol}) can be written in the form
$$r(t,x) = \int_0^ta(s,x)dB^x(s) + \int_0^tB^x(s)\left[\pd ax(s,x)+c(s,x)\right]ds
+ r_0(t+x),$$
where the integral $\int_0^ta(s,x)dB^x(s)$  is the usual Wiener-Ito integral.
\label{wnth}
\end{theorem}
\section{Application to Yield Curves}

Let $r(0,x)=r_0(x)$ be a given yield curve at time $t=0$, that is
the interest on the investment maturing at $x$. Let
$r(t,x)$, $t,x\geq 0$ denote similar curve at time $t$, that is the interest rate at
time $t$ on the investment maturing at $t+x$. In the absence of noise it should hold that
$r(t,x)=r(0,t+x)=r_0(t+x)$, otherwise one can make a riskless profit. Assuming that $r_0(x)$ is smooth,
the  evolution of
the yield curve therefore is described by pde
\begin{equation}
\pd rt-\pd rx=0,
\label{nspde}
\end{equation}
with the initial condition $r(0,x)=r_0(x)$, see Musiela and Sondermann 1994.
Consider now a stochastic analogue of this pde given by (\ref{gspde})
$
\pd rt-\pd rx = DW,
$
where $W(t,x)=B(t,t+x)$.
The reason we take $B(t,t+x)$
rather than  $B(t,x)$ is because we model the yield at the   point $t+x$, and it would be natural
 to write $r(t,t+x)$ instead of $r(t,x)$ and  add to it the
basic noise at the same point
$(t,t+x)$.  The process $W(t,x)=B(t,t+x)=B^x(t)$ for any fixed
$x>0$,   is not a Brownian motion, however it is a continuous martingale
with independent but non-stationary increments.
Condition $W(0,x)=0$ for any $x>0$ in Theorem 2 means here that  at
time zero there is no uncertainty in the yield curve.

 Musiela and Sondermann 1994 introduced the noise by considering for all $x>0$ the stochastic  differential
equations
\begin{equation}
dr(t,x)=\alpha(t,x)dt +\sigma(t,x) dW(t),
\label{MSmodel}
\end{equation}
which hold for all $x\geq 0$.
In comparison  with  this equation,  the
general solution (\ref{gsol})  for any fixed $x>0$ satisfies
$$dr(t,x)=a(t,x)dB^x(t)+\left[B^x(t)\left(\pd ax(t,x)+c(t,x)\right)+r_0'(t+x)\right]dt.$$
Hence in the spde solution the ``drift" term includes random noise, moreover  the equations are driven
 by martingales $B^x(t)$ dependent on $x$.

\section{Proofs}
The proof of Theorem \ref{genPDEth} requires the following basic lemma.
\begin{lemma}
Let $g$ be an integrable function on $[0,+\infty)\times[0,+\infty)$, and
assume that for any test function, on $(0,+\infty)\times(0,+\infty)$, $\varphi$,
\begin{equation}
\dint g(t,x)\pdd\varphi tx(t,x)dx = 0.
\label{lemmahyp}
\end{equation}
Then,  $g(t,x) = g(t,0)+g(0,x)+\mbox{\rm const}$.
\label{lemma}
\end{lemma}

\subsubsection*{Proof of Theorem \ref{genPDEth}}
For any test function $\varphi$,
$-\left<U,\pd\varphi t\right> = <W,D^\star\varphi>,$
or equivalently
\begin{eqnarray*}
\dint U(t,x)\pd\varphi t(t,x)dtdx & = & \dint
W(t,x)\left[\pd{(a\varphi)}t+\pd{(b\varphi)}x-c\varphi\right](t,x)dtdx\\
& = & \dint W(t,x)\left[a\pd\varphi t+b\pd\varphi x+\left(\pd at
+\pd bx-c\right)\varphi\right](t,x)dtdx.
\end{eqnarray*}
Integrating by parts, we get
$$\dint\left\{\int_0^xU(t,y)dy-\int_0^xa(t,y)W(t,y)dy-\int_0^tb(s,x)W(s,x)ds+
\right.$$
$$\left.\int_0^x\left[\int_0^t\left(\pd at+\pd
bx-c\right)(s,y)W(s,y)ds\right]dy\right\}
\pdd\varphi tx(t,x)dtdx = 0.$$
(\ref{solgenPDE}) now follows from Lemma \ref{lemma}.
\hfill$\Box$

\subsubsection*{Proof of Theorem \ref{pspdeth}}
Make a change of variables $\tau= t$, $\xi = t  + x$,
$\rho(\tau,\xi) = r(t,x)$,
$\alpha(\tau,\xi) = a(t,x)$, $\beta(\tau,\xi) = b(t,x)$, $\gamma(\tau,\xi) = c(t,x)$,
$\psi(\tau,\xi) = \varphi(t,x)$ and
$V(\tau,\xi) = W(t,x)$.   Perform the  above change of variables in
(\ref{defgspde}) and use the relation $b(t,x)=-a(t,x)$ to obtain
\begin{equation}
\left<\rho,\pd\psi\tau\right> =
\left<V,\alpha\pd\psi\tau+\left(\pd\alpha\tau-\gamma\right)\psi\right>.
\label{chgvargspde}
\end{equation}
It follows from Theorem \ref{genPDEth} that
$$\int_0^\xi[\rho(\tau,\eta)-\rho(0,\eta)]d\eta
= \int_0^\xi\alpha(\tau,\eta)V(\tau,\eta)d\eta
 +\int_0^\xi\int_0^\tau\left(\pd\alpha t-c\right)(\theta,\eta)V(\theta,\eta)d\theta d\eta.$$
Differentiating with respect to $\xi$, we get
\begin{equation}
\rho(\tau,\xi) =
\alpha(\tau,\xi)V(\tau,\xi) - \int_0^\tau V(\theta,\xi)
\left(\pd\alpha\tau-\gamma\right)(\theta,\xi)d\theta + \rho(0,\xi).
\label{chgvarsol}
\end{equation}
Going back to the original variables and taking
into account the initial condition  we obtain (\ref{gsol}).  Notice that a
formal derivation  of (\ref{chgvarsol}) is obtained by integrating  with
respect to $\tau$ the equality $\pd\rho\tau=\alpha \pd V\tau+\gamma V$
(to be understood  in the sense of generalized functions).
\hfill$\Box$

\subsubsection*{Proof of Theorem \ref{wnth}}
The proof of Theorem \ref{wnth} is based on the following lemmas.

\begin{lemma}
Let $\X$ be a Gaussian measure on a measurable space $(E,{\cal E})$ with
intensity $\mu$, $R$ and $S$ be two functions in $L^2(E,{\cal E},\mu)$, $F$
and $G$ be two measurable sets, and $(F^n_k)_k$ and
$\{G^n_k)_k$ be measurable partitions of $F$ and $G$ respectively.
Assume that
\ben
\item $\mu(F)<+\infty$ and $\mu(G)<+\infty$;
\item $\lim_{n\rightarrow\infty}\sup_{1\leq k\leq n}\mu(F^n_k) =
\lim_{n\rightarrow\infty}\sup_{1\leq k\leq n}\mu(G^n_k) = 0$.
\een
If $(F^n_k\cap G^n_k)_k$ is a partition of $F\cap G$, then
$$\lim_{n\rightarrow\infty}\sum_{k=1}^nR^n_kS^n_k\X(F^n_k)\X(G^n_k) =
\int_{F\cap G}RSd\mu$$
in the $L^2$-sense, where $(R^n_k)$ and $(S^n_k)$ are
approximating sequences on $(F^n_k\cap G^n_k)_k$ of $R$ and $S$ respectively.\\
If, on the other hand, for all $k$ and $l$, one or both of $F^n_k\cap G^n_l$
and  $F^n_l\cap G^n_k$ are empty (in particular $F^n_k\cap G^n_k = \emptyset$),
then
$$\lim_{n\rightarrow\infty}\sum_{k=1}^nR^n_kS^n_k\X(F^n_k)\X(G^n_k) = 0$$
in the $L^2$-sense.
\label{quadlem}
\end{lemma}

\subsubsection*{Proof of Lemma \ref{quadlem}:}
Let us first note that $\var(\X(F^n_k)\X(G^n_k)) = \mu(F^n_k)\mu(G^n_k)
+\mu(F^n_k\cap G^n_k)^2$, and that, for $k\neq l$,
$\cov(\X(F^n_k)\X(G^n_k),\X(F^n_l)\X(G^n_l)) =
\mu(F^n_k\cap G^n_l)\mu(F^n_l\cap G^n_k)$. It follows that, if
$(F^n_k\cap G^n_k)_k$ is a partition of $F\cap G$, then
the random variables $\X(F^n_k)\X(G^n_k)$ and $\X(F^n_l)\X(G^n_l)$,
for $k\neq l$, are independent.
\begin{eqnarray*}
\lefteqn{\norm{\sum_{k=1}^nR^n_kS^n_k\X(F^n_k)\X(G^n_k) -
\int_{F\cap G}RSd\mu}_2\leq}\\
&&\norm{\sum_{k=1}^nR^n_kS^n_k\X(F^n_k)\X(G^n_k) -
\sum_{k=1}^nR^n_kS^n_k\mu(F^n_k\cap G^n_k)}_2 +
\left|\sum_{k=1}^nR^n_kS^n_k\mu(F^n_k\cap G^n_k) - \int_{F\cap G}RSd\mu\right|.
\end{eqnarray*}
The second term of the right hand side goes to zero. By the independence of
the random variables $\X(F^n_k)\X(G^n_k)$,the first term becomes
\begin{eqnarray*}
\var\left(\sum_{k=1}^nR^n_kS^n_k\X(F^n_k)\X(G^n_k)\right)
& = & \sum_{k=1}^n(R^n_kS^n_k)^2\var\left(\X(F^n_k)\X(G^n_k)\right)\\
& = & \sum_{k=1}^n(R^n_kS^n_k)^2[\mu(F^n_k)\mu(G^n_k)+\mu(F^n_k\cap G^n_k)^2]\\
& \leq & 2\sum_{k=1}^n(R^n_kS^n_k)^2\mu(F^n_k\cup G^n_k)^2\\
& \leq & 2\sup_{1\leq k\leq n}\mu(F^n_k\cup G^n_k)
\sum_{k=1}^n(R^n_kS^n_k)^2\mu(F^n_k\cup G^n_k),
\end{eqnarray*}
which completes the proof of the first statement. The second statement follows
in the same way.
\hfill$\Box$

\begin{lemma}
Let $\X$ be a Gaussian measure on a measurable space $(E,{\cal E})$ with
intensity $\mu$, $F$ be a measurable set, and $(F^n_k)_k$
be a measurable partition of $F$. If there exists $\kappa>0$ such that
$$\lim_{n\rightarrow\infty}n^\kappa\sup_{1\leq k\leq n}\mu(F^n_k) = 0,$$
then, with probability one,
$$\lim_{n\rightarrow\infty}\sup_{1\leq k\leq n}|\X(F^n_k)| = 0.$$
\label{ascv}
\end{lemma}

\noindent
{\bf Proof of Lemma \ref{ascv}:} Fix $\varepsilon>0$.
\begin{eqnarray*}
\P\left[\sup_{1\leq k\leq n}|\X(F^n_k)|>\varepsilon\right] & = &
1-\P\left[\sup_{1\leq k\leq n}|\X(F^n_k)|\leq\varepsilon\right]\\
& = & 1-\prod_{k=1}^n\left[2\Phi\left(\frac{\varepsilon}
{\sqrt{\mu(F^n_k)}}\right)-1\right]\\
& \leq & 1-\prod_{k=1}^n\left[1-\exp\left\{-\sqrt{\frac{2}{\pi}}
\frac{\varepsilon}{\sqrt{\mu(F^n_k)}}\right\}\right]\\
& \leq & 1-\left[1-\exp\left\{-\sqrt{\frac{2}{\pi}}
\frac{\varepsilon}{\sqrt{\mu(F^n)}}\right\}\right]^n,
\end{eqnarray*}
where $\mu(F^n) = \sup_{1\leq k\leq n}\mu(F^n_k)$. The proof is ended
by the fact that
$$\sum_n\left(1-\left[1-\exp\left\{-\sqrt{\frac{2}{\pi}}
\frac{\varepsilon}{\sqrt{\mu(F^n)}}\right\}\right]^n\right)<+\infty.$$
\hfill$\Box$

\noindent
{\bf Proof of Theorem \ref{wnth}:} The `` if " part is given in Theorem
\ref{pspdeth}. We now prove the `` only if " part.
Using the notations of the proof of Theorem \ref{pspdeth}, we obtain,
instead of (\ref{chgvargspde}) and without the assumption $b(t,x) = -a(t,x)$,
$$\left<\rho,\pd\psi\tau\right> =
\left<V,\alpha\pd\psi\tau+(\alpha+\beta)\pd\psi\xi+
\left(\pd\alpha\tau+\pd{(\alpha+\beta)}\xi-\gamma\right)\psi\right>.$$
It now follows from Theorem \ref{genPDEth} that for a solution to exist in the space
of functions,
$$Z(t,x) := \int_0^tA(s,x)B(s,s+x)ds =
\int_0^tA(s,x)W(s,x)ds =
\int_0^\tau(\alpha+\beta)(\theta,\xi)V(\theta,\xi)d\theta,$$
where $A=a+b$, must be differentiable in $x$. Let $Y(t,x) = A(t,x)B(t,t+x)$,
$x<y$, $\delta^n = (y-x)/n$ and $\delta^n_k = k\delta^n$. We have,
\begin{eqnarray*}
\lefteqn{\sum_{k=1}^n\left(Z(t,x+\delta^n_k)-Z(t,x+\delta^n_{k-1}\right)^2}\\
& = & \sum_{k=1}^n\left(\int_0^t\left(Y(s,x+\delta^n_k)
- Y(s,x+\delta^n_{k-1})\right)ds\right)^2\\
& = & \int_0^t\int_0^t\sum_{k=1}^n\left[\left(Y(r,x+\delta^n_k)-
Y(r,x+\delta^n_{k-1})\right)
\left(Y(s,x+\delta^n_k)-Y(s,x+\delta^n_{k-1})\right)\right]drds.
\end{eqnarray*}
Now $\bds \Delta Y^n_k(s) := Y(s,x+\delta^n_k)-Y(s,x+\delta^n_{k-1}) =
A^n_k(s)\B(F^n_k(s))+B^n_{k-1}(s)\A(G^n_k(s))\eds$, with $A^n_k(s) =
A(s,x+\delta^n_k)$, $B^n_{k-1}(s)=B(s,s+x+\delta^n_k)$, $F^n_k(s) =
[0,s]\times[s+x+\delta^n_{k-1},s+x+\delta^n_k]$, $G^n_k(s) =
[0,s]\times[x+\delta^n_{k-1},x+\delta^n_k]$, $\B$ is the
Gaussian measure associated to the Brownian sheet $B$, and $\A$ the measure
associated to the $C^1$ function $A$. Thus
\begin{eqnarray*}
\lefteqn{\sum_{k=1}^n\Delta Y^n_k(r)\Delta Y^n_k(s)}\\
& = &\sum_{k=1}^nA^n_k(r)A^n_k(s)\B(F^n_k(r))\B(F^n_k(s))
+\sum_{k=1}^nA^n_k(r)B^n_{k-1}(s)\B(F^n_k(r))\A(G^n_k(s))\\
&&+\sum_{k=1}^nB^n_{k-1}(r)A^n_k(s)\A(G^n_k(r))\B(F^n_k(s))
+\sum_{k=1}^nB^n_{k-1}(r)B^n_{k-1}(s)\A(G^n_k(r))\A(G^n_k(s)).
\end{eqnarray*}
Now each one of the last three terms converges, with probability one, to 0.
Indeed
$$\left|\sum_{k=1}^nA^n_k(r)B^n_{k-1}(s)\B(F^n_k(r))\A(G^n_k(s))\right|
\leq \sup_{1\leq k\leq n}|\B(F^n_k(r))|
\left|\sum_{k=1}^nA^n_k(r)B^n_{k-1}(s)\A(G^n_k(s))\right|,$$
which, according to Lemma \ref{ascv}, goes to 0. Recall that, if $\mu$ is the
intensity of the Gaussian measure $\B$, then
$$\sup_{1\leq k\leq n}\mu(F^n_k(r)) = r\frac{y-x}{n}.$$
The same goes for the third term. The convergence to 0 of the fourth term
follows from
$$\lim_{n\rightarrow\infty}\sup_{1\leq k\leq n}|\A(F^n_k(r))| = 0.$$
To investigate the
asymptotic behavior of the first term, we first make the following observation.
For $r\neq s$, there is $n_0$ $(n_0\geq(y-x)/|r-s|)$ such that for any $n\geq
n_0$, and any $k$ and $l$, $F^n_k(r)\cap F^n_l(s) = \emptyset$ or
$F^n_l(r)\cap F^n_k(s) = \emptyset$. Applying Lemma \ref{quadlem}, we see
that $\sum_{k=1}^nA^n_k(r)A^n_k(s)\B(F^n_k(r))\B(F^n_k(s))$ converges in $L^2$
to 0. On the other hand, if $r=s$, then, by Lemma \ref{quadlem},
$$\lim_{n\rightarrow\infty}
\sum_{k=1}^nA^n_k(s)^2\B(F^n_k(s))^2 = \int_x^yA(s,z)^2sdz$$
in the $L^2$ sense. It follows that
$$\lim_{n\rightarrow\infty}\sum_{k=1}^n\Delta Y^n_k(r)\Delta Y^n_k(s) =
\left\{\begin{array}{lcl}
\int_x^yA(s,z)^2sdz&\mbox{if}&r=s\\
0&\mbox{if}&r\neq s
\end{array}\right.$$
in probability, and that
\begin{equation}
\lim_{n\rightarrow\infty}\sum_{k=1}^n
\left(Z(t,x+\delta^n_k)-Z(t,x+\delta^n_{k-1}\right)^2 =
\int_0^t\left(\int_x^yA(s,z)^2sdz\right)ds
\label{quadZ}
\end{equation}
in probability.
Now assume that $Z(t,x)$, as a process in
$x$, is locally H\"older of order
$\varepsilon>1/2$, then
$$\left|Z(t,x+\delta^n_k)-Z(t,x+\delta^n_{k-1})\right| \leq \mbox{const }
(\delta^n)^\varepsilon,$$
and
$$\sum_{k=1}^n\left(Z(t,x+k\delta_n)-Z(t,x+(k-1)\delta_n)\right)^2  \leq
\mbox{const}^2\ (y-x)(\delta^n)^{2\varepsilon-1}\longrightarrow0,\ \ \mbox{as }
n\rightarrow\infty.$$
This combined with (\ref{quadZ}) shows that $A(s,x)$ must be nil.
\hfill$\Box$

\section*{References}

\begin{enumerate}

\item Carmona, R. and Nualart, D. (1988) ``Random Non-linear Wave equations:
Smoothness of the solutions." {\it Probab. Th.  Rel. Fields} 79, 469-508.

\item Freidlin, M.I. (1988) ``Random  Perturbations of Reaction-Diffusion Equations: The
Quasy-Deterministic Approximation," {\it Trans. Amer. Math. Soc. } 305, 665-697.





\item Musiela, M.  and  Sondermann, D. (1994)  ``Different Dynamical
Specifications of the Term Structure of Interest Rates and Their Implications,"
{\it Workshop on Stochastics and Finance, ANU, Canberra, Statistics Research
Report No SRR 016-94}.

\item Lang, S. {\it Real and Functional Analysis}, Springer 1993.

\item Revuz, D. and Yor, M. {\it Continuous Martingales and Brownian Motion},
Springer 2001.
\item Rozanov, Y. A. {\it Random fields and stochastic partial differential equations}, Nauka 1995, Kluwer, 1998.

\item Walsh, J.B. ``An Introduction to Stochastic Partial Differential
Equations," {\it Lecture Notes in Math., } vol. 1180, 265-437,  Springer, 1986.
\end{enumerate}


\end{document}